\documentclass{article}

\advance\textheight -\headheight
\advance\textheight -\headsep
\oddsidemargin\paperwidth
\advance\oddsidemargin -\textwidth
\divide\oddsidemargin2 
\ifdim\oddsidemargin<.5truein \oddsidemargin.5truein \fi
\advance\oddsidemargin -1truein
\evensidemargin\oddsidemargin
\topmargin\paperheight \advance\topmargin -\textheight
\advance\topmargin -\headheight \advance\topmargin -\headsep
\divide\topmargin2
\ifdim\topmargin<.5truein \topmargin.5truein \fi
\advance\topmargin -1.4truein\relax

\usepackage{amsfonts, amsmath, amsthm, amssymb}
\usepackage{graphicx}

\usepackage{hyperref}
\hypersetup{
    colorlinks=true,
    linkcolor=blue,
    citecolor=blue,
    urlcolor=blue
}

\newtheorem{theorem}{Theorem}

\newtheorem{corollary}[theorem]{Corollary}

\numberwithin{equation}{section}

\begin{document}

\title{Fullerenes and Belyi functions}

\author{
Nikolai M. Adrianov\footnote{
Lomonosov Moscow State University, 119991 Moscow, Russia, e-mail: {\tt nadrianov@gmail.com}},
George B. Shabat\footnote{
Lomonosov Moscow State University, 119991 Moscow, Russia;
Russian State University for the Humanities, 125993 Moscow, Russia, e-mail: {\tt george.shabat@gmail.com}}
}

\date{}

\providecommand{\keywords}[1]{\textbf{Keywords:} #1}

\maketitle

\begin{abstract}
The paper is an attempt to apply the theory of dessins d'enfants to the theory of fullerenes. The classical results concerning the calculation of the dodecahedron Belyi function are presented and then applied to the calculation of the Belyi function of the barrel, and the euclidean geometry of the latter is investigated. The non-existence of the fullerene with the only hexagonal face is established by the methods of dessins d'enfants.
\end{abstract}

\tableofcontents

\setcounter{section}{-1}

\section{Introduction}

Fullerenes, being a hot point of chemistry and other natural sciences, have lots of  applications, ranging from quantum computers (see, e.g., \cite{BenjaminETAL2006}) to medicine (e.g., \cite{PetrovicETAL2015}).

The mathematical aspects of the fullerenes theory are also diverse. Various fullerenes-related topics of graph theory, finite groups and their representations, geometry of polyhedra, isoperimetric problems, etc., can be found (among many other papers) in \cite{ChuSte1993} and \cite{DeSiSh2013}. The present paper is aimed at introducing yet another mathematical discipline that, hopefully, can promote the understanding of combinatorics and geometry of fullerenes -- the {\it dessins d'enfants} theory; this theory was created by Alexander Grothendieck in \cite{Gro1984}, and the modern introductions can be found in \cite{LanZvo2004} and \cite{GG2012}.

The main ingredient that can be contributed to fullerenes theory by the dessins d'enfant theory is the {\it global conformal structure} of the fullerenes graphs on the Riemann sphere; more precisely, we have a 3-parametrric real family of {\it metric} structures. One of the current central questions is -- whether these structures have some relations to the ones observable in Nature. In some cases, to one of which (the {\it barrel}) a section of this paper is devoted, a certain {\it distinguished} metric  structure can be uniquely selected from the above family. If the answer to the central question turns out to be affirmative, the new field of the {\it applied dessins d'enfants theory} will be opened; some initial problems of this field are listed in the conclusion of the present paper.

Hopefully, the approach we develop has some relations with the rich  mathematics around fullerenes  which we have learned about from the recent E. Katz article \cite{Katz-2014}. It turned out that some fullerene-related chemical structures ``solve'' certain non-trivial combinatorial and variational problems, possibly related, e.g., to the {\it dense circle packing} (see \cite{BoSte2004}) which is one of  the methods of approximate solution of the problem we are dealing with in the present paper. E. Katz refers to the monograph \cite{MB-1936} in which its author, Dmitry Morduhai-Boltovskoi\footnote{This mathematician with rather diverse interests is widely known for his translation of Euclid's Elements. Also the writer A. Solzhenitsyn was his student.}, applied various mathematical tools to describe the geometry of {\it radiolaria} discovered in 19-th century (and studied deeply by Alan Turing and his students) and to the fullerene-like structures, decades before the genuine fullerenes were discovered.

In the section 1 we briefly introduce the elements of dessins d'enfants theory in their relation to spherical graphs and, in particular, the fullerenes; the corresponding Belyi functions are introduced. The section 2 contains the main algebraic equation solving which provides these Belyi functions. 

In the section 3 we discuss the  fullerenes with the smallest numbers of  hexagons. As for the one without hexagons, it is just the {\it dodecahedron}. The calculation of its Belyi function (of course, in different terminology) belong to the 19th century mathematics. However, we perform the detailed calculation in terms of dessins d'enfant; it turns out to be much shorter (our several pages vs. the considerable part of the famous Klein's book) and, perhaps, conceptually clearer. Then we prove the non-existence of a fullerene with only one hexagon. It is the well-known folklore result; the overview of the recent papers containing the rigorous proof of this fact is given in the subsection~3.1. We also demonstrate that dessins provide a clear -- though somewhat cumbersome -- algebraic way of dealing with the problems of this kind. (While working on this proof we found out that the similar problems were considered around the edge of 19th and 20th centuries by several French and Russian mathematicians, some of them currently little-known; we give the corresponding formulations and references).  Then we study in some detail the fullerene with two hexagons, or the so-called {\it barrel}. Using some symmetry arguments we mostly reduce it to the dodecahedron calculations. In the end of the section we study the metric of the barrel with the hope that it can be compared with the barrels that can be found in Nature. 

In the short final section 4 we speculate about the possible relations between the dessins d'enfants theory and chemistry.

We are indebted to Victor Buchshtaber, whose lecture in the Dubna summer school (2015) and subsequent discussions involved us in this exciting domain.

\section{Fullerenes and dessins d'enfants}

We start with establishing relations between the traditional objects of the mathematical theory of fullerenes and those of dessins d'enfants theory.

{\bf 1.0. Brief overview of dessins d'enfants.} 
{\it Dessin d'enfant} is a bicolored graph $\Gamma$ embedded into a compact connected oriented surface $X$ in such a way that $X\setminus\Gamma$ is homeomorphic to a disjoint union of open discs. These discs are called the {\it faces} of the dessin.

The term was coined by Alexander Grothendieck, these objects are known since 19th century as {\it bicolored maps} (or {\it hypermaps}) {\it on surfaces}. Dessins d'enfants theory differs from the others, considering the same objects, by the {\it categorical approach} and establishing equivalences between the appropriately defined category $\mathcal{DESSINS}$ and seemingly very different categories. For the needs of the present paper we mention just one of them -- the category $\mathcal{BELPAIR}(\mathbb{C})$ of {\it Belyi pairs} over complex numbers. The objects of the latter are pairs $(\mathfrak{X},\beta)$, where $\mathfrak{X}$ is a smooth compact connected Riemann surface and $\beta$ a non-constant meromorphic function on it whose finite critical values are only 0 and 1.

To a Belyi pair a dessin d'enfant is associated in the following manner:
the surface $X:=\mathbf{top}(\mathfrak{X})$ is defined as the {\it topological model} of the Riemann surface (that simply means forgetting the complex structure) and the graph $\Gamma:=\beta^{-1\circ}([0,1])$ is defined as the pre-image of the segment connecting the finite critical values; black and white vertexes of the graph $\Gamma$ are the pre-images of $0$ and $1$ respectively.

We have given a sketch of the definition of the functor
$$
\mathbf{draw}: \mathcal{BELPAIR}(\mathbb{C})\longrightarrow\mathcal{DESSINS}.
$$
It turns out to be {\it category equivalence}; the details can be found in \cite{ShaVoe1990} and many other places. This means, in particular, that \underline{any} dessin d'enfant corresponds to a certain complex-analytical structure. By the {\it Riemann existence theorem} this complex structure has the algebro-geometric nature: $\mathfrak{X}$ is an algebraic curve and $\beta$ a rational function on it.

In the present paper we restrict ourselves by the surfaces $X\simeq\mathbf{S}^2$ of genus zero.
According to the above, we identify the 2-sphere with the projective line (Riemann sphere) $\mathbf{P}_1(\mathbb{C})$, and our problem reduces to finding such a rational function $\beta(z)$ that the full pre-image $\beta^{-1\circ}[0,1]$ is isotopically equivalent to $\Gamma$. It is ``well-known''\footnote{This statement is the so-called ``easy Belyi theorem''; however, the experts are of different opinions concerning its obviousness.} \cite{LanZvo2004} that one can find the desired function $\beta\in\overline{\mathbb{Q}}(z)$ with {\it algebraic} coefficients. A wide variety of advanced methods starting with \cite{Cou1994} have been developed to calculate Belyi functions; an excellent survey of these techniques can be found in \cite{SijVoi2013}. A fascinating method based on the above-mentioned dense disks packing was developed in \cite{BoSte2004}.

The ramification of a Belyi function $\beta$ over three critical values $0$, $1$ and $\infty$ defines three partitions of $n=\deg \beta$. We call the triple of these partitions a {\it passport} of the dessin. Combinatorially the passport is defined by the degrees of black vertexes, white vertexes and faces of the dessin. See \cite{LanZvo2004} for more details.

{\bf 1.1. Polyhedra and spherical graphs.}
For any convex polyhedron ${\mathcal P}$ we can consider the projection of the vertexes and edges of ${\mathcal P}$ from an internal point $O$ of ${\mathcal P}$ onto any sphere centered at $O$, thus getting a connected graph on the sphere; a stereographic projection of this graph gives us a plane graph called {\it Schlegel diagram}. The famous {\it Steinitz theorem} \cite{Ste1922} states that any simple\footnote{Simple graph is a graph without loops and multiple edges.} planar 3-connected graph is the Schelegel diagram of a polyhedron. A compelling insight into the history of the Steinitz theorem can be found in \cite{Grunbaum-2007}.

\begin{figure}[ht]
\begin{center}
\includegraphics[scale=1.0]{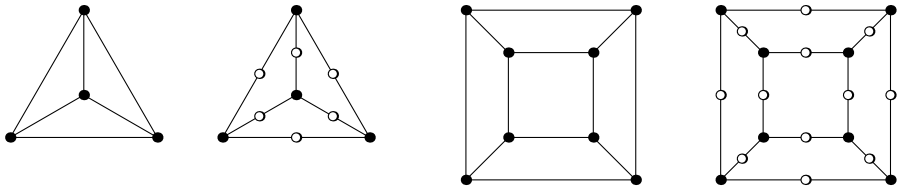}
\end{center}
\caption{Schlegel diagrams and the corresponding dessins.}
\label{fig:Schlegel}
\end{figure}

Schlegel diagram is not necessarily bicolorable, so we put a white vertex in the middle of every edge; the result is a dessin with white vertexes of degree 2, see fig.~\ref{fig:Schlegel}. It follows from 3-connectedness of the graph that the degrees of black vertexes are $\ge 3$. The passport of the tetrahedron dessin is $(3^4\,|\,2^6\,|\,3^4)$ and the passport of the cube dessin is $(3^8\,|\,2^{12}\,|\,4^6)$.

Schlegel diagrams provide a convenient language to speak about the combinatorial structure of polyhedra as soon as we are not interested in their metric properties (see e.g. \cite{BuchEro2015}, \cite{BuchEro2017}). Nevertheless we would like to point out that any dessin defines a unique complex structure on the sphere.

Any connected plane 3-valent graph all of whose faces are 5- and 6-gons can be realized by a polyhedra (fullerene) (see \cite{BuchEro2017}, \cite{DeSiSh2013}). Therefore dessins d'enfants with passports $(3^{2n}\,|\,2^{3n}\,|\,5^{12}\, 6^{n-10})$ correspond exactly to the fullerenes with $2n$ vertexes. Note that without the 5, 6-assumption on the faces that is not true: not every realization of a passport by a dessin corresponds to a polyhedron.

{\bf 1.2. Face vectors.} We use the traditional notations $f_i$ for $i=0,1,2$ the number of faces of dimension $i$ of the polyhedron corresponding to a fullerene. Then the following system of equations holds:
$$
\left\{
\begin{aligned}
  f_0-f_1+f_2&=2 && (\text{\it Euler})\\
  3f_0&=2f_1&& (\text{\it trivalent})\\
  f_2&=p_5+p_6&& (\text{\it face control})\\
  3f_0&=5p_5+6p_6&& (\text{\it Fuller})\\
\end{aligned}
\right.
$$
The formal operations with this systems imply the well-known fullerene relation
$$
\boxed{ p_5=12}
$$
and
$$
\left\{
\begin{aligned}
  f_0&=20+2p_6&& (vertexes)\\
  f_1&=30+3p_6&& (edges)\\
  f_2&=12+p_6 && (faces)\\
\end{aligned}
\right.
$$
These relations are easily memorisable: the case $p_6=0$ corresponds to the dodecahedron. Gr\"unbaum and Motzkin \cite{GrunbaumMotzkin-1963} proved that all the values $p_6\ge 2$ are realizable.

\section{Belyi functions corresponding to fullerenes}

According to subsection {\bf 1.0}, the problem of finding the Belyi pair, corresponding to a fullerene with $p_6$ hexagons reduces to finding the rational Belyi function $\beta\in\mathbb{C}(z)$ of the form
\begin{equation}
\label{eq:2a}
\beta=k\frac{V^3}{P^5H^6},
\end{equation}
where $k\in\mathbb{C}\setminus\{0\}$ and $V,P, H\in\mathbb{C}[z]$ are the {\it vertex}, {\it pentagon} and {\it hexagons} polynomials.

According to subsection {\bf 1.2},
\begin{equation}
\label{eq:2b}
\deg P=12,\quad \deg H=p_6,
\end{equation}
\begin{equation}
\label{eq:2c}
\deg V=20+2p_6.
\end{equation}
The function $\beta$ being Belyi means that $\beta-1$ has all the zeros of multiplicity exactly 2, hence, calculating
$$
\beta-1=k\frac{V^3}{P^5H^6}-1=\frac{kV^3-P^5H^6}{P^5H^6},
$$
we arrive at the main equation
\begin{equation}
\label{eq:2d}
\boxed{kV^3-P^5H^6=M^2}
\end{equation}
where $M$ is a {\it midpoint} polynomial, assumed to have no multiple roots.

The unknowns  of the main  equation (\ref{eq:2d}) are
$$
k;\text{ coefficients of } V, P, H, M.
$$
Assuming  $V$, $P$ and $H$ to be monic, count the unknowns:
\begin{equation}
\label{eq:2e}
1 + (20+2p_6) + 12 + p_6 + (31+3p_6) = 64+6p_6.
\end{equation}
The number of equations is
\begin{equation}
\label{eq:2f}
1+2\deg M=1+2(30+3p_6)=61+6p_6.
\end{equation}
The excess by (\ref{eq:2e}) and (\ref{eq:2f}) equals
$$
64-61=3=\dim\mathsf{PSL}_2(\mathbb{C}).
$$
It corresponds to the freedom of fractional-linear transformations of the variable $z$.

\section{Smallest fullerenes and their Belyi functions}

\setcounter{subsection}{-1}
\subsection{$p_6=0$. Dodecahedron}

\begin{figure}[ht]
\begin{center}
\includegraphics[scale=1.0]{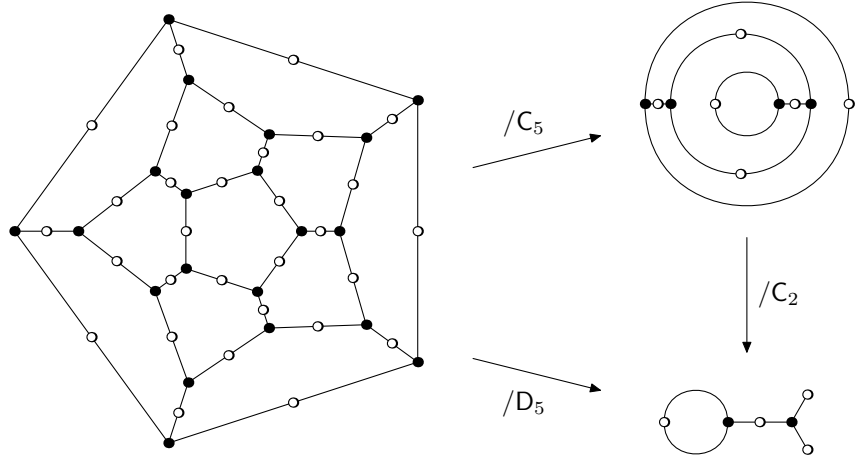}
\end{center}
\caption{Dodecahedron dessin and its factors by $\mathsf{C}_5$ and $\mathsf{D}_5$.}
\label{fig:dodecahedron}
\end{figure}

The smallest fullerene $C_{20}$ with $p_6=0$ is one of the Platonian solids, dodecahedron. It has
$$
\boxed{\#\text{vertexes}=20,\ \#\text{edges}=30,\ \#\text{faces}=12}
$$
which is consistent with the formulas from subsection {\bf 1.2}. The corresponding dessin $D_{60}$ is shown at the left in fig.~\ref{fig:dodecahedron}. In this subsection we compute its Belyi function. 

The Belyi function corresponding to the dodecahdron
was actually found in a classical work of Schwarz \cite{Schwarz1873} (in 1873, exactly 150 years ago!). The same formulas can be found in the Klein's book \cite{Klein1888} published in 1884. We believe that the pictorial approach suggested by the theory of dessins d'enfants can clarify these calculations.

\begin{figure}[ht]
\begin{center}
\includegraphics[scale=1.0]{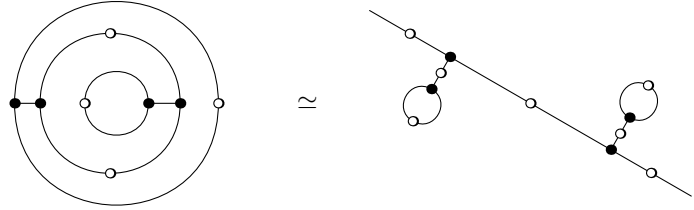}
\end{center}
\caption{Two images of the dessin $D_{12}$.}
\label{fig:12edges-factor}
\end{figure}

The automorphism group of the dodecahedron is isomorphic to $\mathsf{A}_5$, the alternating group on 5 points. Let us consider the factor of $D_{60}$ by any of the subgroups of $\mathsf{A}_5$ isomorphic to $\mathsf{D}_5$. This factorization can be performed in two steps: first factor the dodecahedron by $\mathsf{C}_5$ which will result in dessin $D_{12}$ with 12 edges and the passport $(3^4\,|\,2^6\,|\,5^2\,1^2)$, see fig.~\ref{fig:dodecahedron}. The dessin $D_{12}$ has a non-trivial automorphism of order 2, which is not so easy to see when we draw it as in fig.~\ref{fig:dodecahedron}; however, it becomes obvious if we draw $D_{12}$ alternatively, see fig.~\ref{fig:12edges-factor}. Factorizing $D_{12}$ by this automorphism we get dessin $D_{6}$ with 6 edges and passport $(3^2\,|\,2^2\,1^2\,|\,5^1\,1^1)$.

Let $\beta_6(z)$ be a Belyi function corresponding to dessin $D_6$. To make calculations simpler we assume that its only pole of degree 5 is placed at $z=\infty$ and the only pole of degree 1 is at $z=0$. This leaves us only the multiplicative freedom $z \to \alpha z$. Then the desired function takes the form
$$
\beta_6(z) = \frac{\left(z^2+a_1z+a_0\right)^3}{kz} =
\frac{\left(z^2+b_1z+b_0\right)^2 \left(z^2+c_1z+c_0\right)}{kz} + 1
$$
and the following polynomial is equal to zero:
$$
S := \left(z^2+a_1z+a_0\right)^3 - \left(z^2+b_1z+b_0\right)^2 \left(z^2+c_1z+c_0\right) - kz = 0,
$$
where $(z^2+a_1z+a_0)$, $(z^2+b_1z+b_0)$ and $(z^2+c_1z+c_0)$ are pairwise co-prime.

Coefficients of the polynomial $S$ at $z^5$, $z^4$ are linear in $c_1$, $c_0$, so we obtain
\begin{equation}
\label{eq:dodecahedron-1}
\begin{array}{lll}
c_1 &=& 3a_1 - 2b_1,\\[2pt]
c_0 &=& 3a_0 - 2b_0 + 3(a_1 - b_1)^2.\\[3pt]
\end{array}
\end{equation}

If the equality $b_1=a_1$ would hold then we would have
$$
S = (a_0 - b_0)^2 (3z^2 + 3 a_1 z + a_0 + 2 b_0) - kz,
$$
and it would imply $b_0=a_0$, so $(z^2+a_1z+a_0) = (z^2+b_1z+b_0)$ would not be co-prime. Therefore $b_1\ne a_1$.

Substituting (\ref{eq:dodecahedron-1}) into the coefficient of $S$ at $z^3$ we get
$$
(a_1 - b_1)(a_1^2 - 5 a_1 b_1 + 4 b_1^2 + 6 a_0 - 6 b_0) = 0
$$
and hence
\begin{equation}
\label{eq:dodecahedron-2}
b_0 =a_0 + \frac{1}{6}(a_1 - b_1)(a_1 - 4 b_1).
\end{equation}
Substituting (\ref{eq:dodecahedron-1}) and (\ref{eq:dodecahedron-2}) into the coefficient 
of $S$ at $z^2$ we get
$$
-\frac{(a_1 - b_1)^2 (11 a_1^2 - 40 a_1 b_1 + 20 b_1^2 + 36 a_0)}{12} = 0
$$
and hence
\begin{equation}
\label{eq:dodecahedron-3}
a_0 = \frac{a_1^2}{4} - \frac{5}{9}(a_1 - b_1)^2.
\end{equation}

Substituting (\ref{eq:dodecahedron-1}), (\ref{eq:dodecahedron-2}) and (\ref{eq:dodecahedron-3}) into the free term of $S$ we get
$$
-\frac{(2 a_1 - 5 b_1)(a_1 - b_1)^5}{27} = 0
$$
and hence
\begin{equation}
\label{eq:dodecahedron-4}
b_1 = \frac{2}{5} a_1.
\end{equation}

Using the multiplicative freedom $z \to \alpha z$ we can set $a_1=10$, then the Belyi function takes the form
\begin{equation}
\label{eq:dodecahedron-5}
\beta_6(z) = \frac{\left(z^2 + 10 z + 5\right)^3}{1728 z} =
\frac{\left(z^2 + 4 z - 1\right)^2 \left(z^2 + 22 z + 125\right)}{1728 z} + 1.
\end{equation}

The roots of the polynomial $(z^2 + 22 z + 125)$ are $-11\pm i$. Let $\mu_1 \in \mathsf{PSL}_2({\mathbb C})$ be such that it maps $0$, $1$ and $\infty$ to $-11+2i$, $0$ and $-11-2i$ respectively:
$$
\mu_1(z) = \frac{-125i(z-1)}{(2+11i)z+(2-11i)}.
$$
Then $\beta_6\big(\mu_1(z^2)\big)$ is a Belyi function of $D_{12}$ such that the 
poles are placed at $\pm 1$. 
Combining it with
$$\
\mu_2(z) = \frac{iz-1}{iz+1}
$$
that maps $0$ to $-1$ and $\infty$ to $1$ we get it in the form where the 
poles of order 1 are placed at $0$ and $\infty$:
$$
\beta_{12}(z) = \beta_6 \circ \mu_1 \circ(z\to z^2) \circ \mu_2 =
\frac{(z^4 + 228 z^3 + 494 z^2 - 228 z + 1)^3}{1728 z (z^2 - 11 z - 1)^5}={}
$$
\begin{equation}
\label{eq:dodecahedron-6}
{}=
\frac{(z^2 + 1)^2(z^4 - 522 z^3 - 10006 z^2 + 522 z + 1)^2}{1728 z (z^2 - 11 z - 1)^5} + 1.
\end{equation}

Now we can write down the Belyi function of the dodecahedron
\begin{equation}
\label{eq:dodecahedron-7}
\beta_{60}(z) = \beta_{12}(z^5) = \frac{(z^{20} + 228 z^{15} + 494 z^{10} - 228 z^5 + 1)^3}
{1728 z^5 (z^{10} - 11 z^5 - 1)^5}. 
\end{equation}

Note that this Belyi function was actually computed by Schwarz while composing the famous {\it Schwarz's list} (of algebraic hypergeometric functions): in his notations (see~\cite{Schwarz1873}, p.330)
$$
\begin{array}{lll}
\varphi_{12} &=& s(1 - 11 s^5 - s^{10}),\\
\varphi_{20} &=& 1 + 228 s^5 + 494 s^{10} - 228 s^{15} + s^{20},\\
\varphi_{30} &=& 1 - 522 s^5 - 10005 s^{10} - 10005 s^{20} + 522 s^{25} + s^{30}
\end{array}
$$
and the Belyi function can be written as
$$
\beta_{60}(-s) = \frac{\varphi_{20}^3}{1728\cdot \varphi_{12}^5} = \frac{\varphi_{30}^2}{1728\cdot \varphi_{12}^5} + 1.
$$

The same function appeared in the Klein's research
devoted to the icosahedron and the solution of equations of the fifth degree, see~\cite{Klein1888}, p.61\footnote{
Russian edition contains interesting modern appendixes by Serre, Arnold, Tyurin. Note however that there is a misprint in this edition: in one of the key formulas (see the $\varphi_{30}$ in the main text) the number 1005 is printed instead of 10005.}.

\subsection{$p_6=1$. Non-existing fullerene}
\label{sec:Non-existing}

It is a common knowledge that there exist no fullerene $C_{22}$ with $p_6=1$. Probably Goldberg \cite{Goldberg-1935} in 1935 was the first who noted that there exist no polyhedra with these parameters. He says: ``It is easy to show that no medial polyhedra exist for $n=11$ and for $n=13$''\footnote{Fullerene with $p_6=1$ is a medial polyhedra with $n=13$ faces in Goldberg's terminology.}. Goldberg makes no attempt to prove it.

Some authors erroneously claim that non-existense of this fullerene was proved by Gr\"unbaum and Motzkin \cite{GrunbaumMotzkin-1963}. That is not correct: Gr\"unbaum and Motzkin prove that there exist a fullerene for any $p_6\ge 2$, but they do not try to prove that there is no one for $p_6=1$. Another proof of the existence is given, e.g., in \cite{VoySte2005}, though the problem of non-existence is still evaded by the authors.

The problem was presented as an {\it Analytical Challenge} puzzle to the readers of {\it Analytical and Bioanalytical Chemistry} \cite{Meija-2006} in 2006 and the challenger gives the following justification \cite{Meija-2006a} \footnote{Also he refers to the paper of Gr\"unbaum and Motzkin as a source of ``more rigorous'' proof.}:

\medskip

\hfill\begin{minipage}{\dimexpr\textwidth-1cm}
{\it Experimentation with fullerene graphs (planar projections) reveals that it is indeed impossible to construct a 22-vertex polyhedron with 12 pentagons and only one hexagon. $\langle...\rangle$ An attempt to construct a fullerene with only one hexagon at the center (Fig. 1 right) fails because the final (encompassing) face of the graph must be another hexagon, thus resulting in a polyhedron with a total of two hexagons and 24 vertexes.}
\end{minipage}

\medskip

\noindent Clearly, it can not be considered as a mathematical proof and this kind of reasoning is difficult to formalize.
A complete, mathematically rigorous, combinatorial proof is given in \cite{BuchEro2015a} (Prop. 1.26).
The current subsection is, to a large extent, inspired by a question of V. Buchstaber, who asked for a proof of non-existence of $C_{22}$ by the methods of the dessin d'enfants theory. In fact, there exist a number of proofs, all of them are based on rather deep mathematical theories.

The main result of this subsection is the following theorem (we do not state it in maximal generality as it goes out of the scope of the paper). In maximal generality it is stated by Vidunas and Filipuk \cite{Vidunas-2014}, Lemma 5.1a. Their proof is short, but assumes that the reader knows the theory of Fuchsian differential equations.
Yet another short proof using the metric structures is given by Izmestiev \cite{Izmestiev-2015}.
We present here a self-contained proof which uses elementary reasonings only (see also the discussion of generalized Fermat equation and a Halphen's result at the end of the subsection). And, finally, one more proof via pure group-theoretical arguments is given in \cite{Adrianov-2023}.

\begin{theorem}
\label{th:Non-realizable}
There exist no dessins of genus $g=0$ with passport $( 3^k \,|\, 2^l \,|\, 5^m\,s^1 )$ for $s\ne 5$.
\end{theorem}

\begin{corollary}
There exist no fullerene $C_{22}$ with $p_6=1$.
\end{corollary}

\proof Suppose there exists a dessin with passport $( 3^k \,|\, 2^l \,|\, 5^m\,s^1 )$. The number of edges of the dessin equals $n=3k=2l=5m+s$. From the Euler's formula we have
$$
k + l + m + 1 - n = 2,
$$
or
$$
\frac{n}{3} + \frac{n}{2} + \frac{n-s}{5} = n + 1,
$$
and it follows that
\begin{equation}
\label{eq:V3M2P5-degree}
\begin{array}{ll}
n = 30 + 6s, & k = 10 + 2s, \\
l = 15+3s, & m = 6 + s.
\end{array}
\end{equation}
Put the center of the only face of degree $s$ to the infinite point of 
$\mathbf{P}_1$ (thus leaving the freedom of affine substitution $z\to az+b$). Then the Belyi function will take the form
\begin{equation}
\label{eq:V3M2P5-Belyi}
\beta=\frac{V^3}{kP^5} = \frac{M^2}{kP^5} + 1,
\end{equation}
where $V$, $M$ and $P$ are monic pairwise coprime polynomials with no multiple roots and
\begin{equation}
\label{eq:V3M2P5-Belyi-deg}
\deg V = k, \quad \deg M = l, \quad \deg P = m.
\end{equation}

Apply the so-called {\it differential trick}, i.e. calculate the derivatives in two ways:
$$
\beta' = \frac{V^2(3V'P-5VP')}{kP^6} = \frac{M(2M'P-5MP')}{kP^6}
$$
and hence
\begin{equation}
\label{eq:V3M2P5-ODE-1}
V^2(3V'P-5VP')=M(2M'P-5MP').
\end{equation}
Simple calculation shows that $3V'P-5VP'$ is a polynomial of degree
$$
k+m-1 = 10 + 2s + 6 + s - 1 = 15 + 3s = \frac{n}{2}
$$
with the leading coefficient equal to $3k-5m = n - (n-s) = s$. Since $V$ and $M$ are coprime,
$M$ divides $3V'P-5VP'$, so both polynomials are of degree $\frac{n}{2}$ and $M$ is monic, so we get
\begin{equation}
\label{eq:V3M2P5-sM}
sM = (3V'P-5VP').
\end{equation}
Similarly
\begin{equation}
\label{eq:V3M2P5-sV2}
sV^2=2M'P-5MP'.
\end{equation}

Substituting $M$ gotten from (\ref{eq:V3M2P5-sM}) into (\ref{eq:V3M2P5-sV2}) we get
\begin{equation}
\label{eq:V3M2P5-ODE-2}
s^2 V^2 = 6 V'' P^2 - 19 V' P' P - 10 V P P'' + 25 V P'^2,
\end{equation}
which is equivalent to
\begin{equation}
\label{eq:V3M2P5-ODE-2a}
V (s^2 V + 10 P P'' - 25 P'^2) = P (6 V'' P - 19 V' P').
\end{equation}

Since $V$ and $P$ are coprime, we conclude that
\begin{equation}
\label{eq:V3M2P5-VR}
VR = 6 V'' P - 19 V' P',
\end{equation}
\begin{equation}
\label{eq:V3M2P5-PR}
PR = s^2 V + 10 P P'' - 25 P'^2
\end{equation}
for some $R\in{\mathbb C}[z]$. From (\ref{eq:V3M2P5-VR}) we infer that
$\deg R \le \deg P - 2$ (it is easy to check that the equality holds, checking the coefficient at $z^{m-2}$, but we do not need this).

From (\ref{eq:V3M2P5-PR})
\begin{equation}
\label{eq:V3M2P5-V-PR}
s^2 V = PR - 10 P P'' + 25 P'^2,
\end{equation}
and substituting this into (\ref{eq:V3M2P5-VR}) we get
$$
\begin{array}{l}
P R^2 - 10 P P'' R + 25 P'^2 R = 6 P P'' R - 7 P P' R' + 6 P^2 R'' + {}\\[6pt]
\qquad{} + 240 P P''^2 + 370 P P' P''' - 60 P^2 P^{IV} - 19 P'^2 R - 760 P'^2 P'',
\end{array}
$$
which is equivalent to
$$
\begin{array}{l}
4 P'^2 (11 R + 190 P'') ={}\\[6pt]
\qquad{}= -P (7 P' R' - 6 P R'' - 370 P' P''' + 60 P P^{IV} + R^2 - 16 P'' R - 240 P''^2).
\end{array}
$$
Since $P$ has no multiple roots, 
$P$ and $P'$ are coprime and hence $P$ has to divide $11 R + 190 P''$, but the degree of the latter is less than $\deg P$. So 
\begin{equation}
\label{eq:V3M2P5-ODE-3}
11 R + 190 P'' = 0
\end{equation}
and
\begin{equation}
\label{eq:V3M2P5-ODE-4}
7 P' R' - 6 P R'' - 370 P' P''' + 60 P P^{IV} + R^2 - 16 P'' R - 240 P''^2 = 0.
\end{equation}
Substituting $11 R = - 190 P''$ into (\ref{eq:V3M2P5-ODE-4}), we get a differential equation of order 4 in $P$:
\begin{equation}
\label{eq:V3M2P5-ODE}
22 P P^{IV} + 45 P''^2 - 66 P' P''' = 0.
\end{equation}
Substituting $11 R = - 190 P''$ into (\ref{eq:V3M2P5-V-PR}) and substituting the result into (\ref{eq:V3M2P5-sM})
we express $V$ and $M$ in terms of $P$:
\begin{equation}
\label{eq:V3M2P5-V}
\frac{11s^2}{25}V = -12 P P'' + 11 P'^2,
\end{equation}
\begin{equation}
\label{eq:V3M2P5-M}
\frac{11s^3}{25}M = 90 P P' P'' - 36 P^2 P''' - 55 P'^3.
\end{equation}

Now we proceed solving the differential equation in polynomial $P$. The left hand side of (\ref{eq:V3M2P5-ODE}) is a polynomial of degree $\le 8+2s$ and the coefficient at $z^{8+2s}$ equals 
$$
\begin{aligned}
22m(m-1)(m-2)(m-3) + 45m^2(m-1)^2 - 66 m^2(m-1)(m-2) ={}\\[2pt]
{}= m(m - 1)(m - 11)(m - 12) = (s-6)(s-5)(s+5)(s+6).
\end{aligned}
$$
Therefore the equation (\ref{eq:V3M2P5-ODE}) has no solutions for $s\ne 5,6$.

Write down $P$ with indeterminate coefficients:
\begin{equation}
\label{eq:V3M2P5-P-indeterminate}
P = z^m + a_{m-1}z^{m-1} + a_{m-2}z^{m-2} + \ldots a_1 z + a_0.
\end{equation}
Using affine substitution $z\to z+b$ (and leaving the multiplicative freedom $z\to az$) we can assume that $a_{m-1}=0$. Now we consider two cases $s=5$ and $s=6$.

{\bf Case $s=5$.} In this case $m=11$. Plugging (\ref{eq:V3M2P5-P-indeterminate}) into (\ref{eq:V3M2P5-ODE}) and solving sequentially linear equations in $a_i$ we get
$$
a_9 = a_8 = a_7 = a_5 = a_4 = a_3 = a_2 = 0,
\quad
a_1 = -\frac{a_6^2}{121},
\quad
a_0 = 0.
$$
Using multiplicative freedom $z\to az$ we can take $a_6=11$ and we get
$$
P = z^{11} - 11 z^6 - z
$$
and from (\ref{eq:V3M2P5-V}) and (\ref{eq:V3M2P5-M})
$$
V = z^{20} + 228 z^{15} + 494 z^{10} - 228 z^5 + 1,
$$
$$
M = z^{30} - 522 z^{25} - 10005 z^{20} - 10005 z^{10} + 522 z^5 + 1.
$$

We have got the formulas for the Belyi function of dodecahedron:
$$
\beta = \frac{(z^{20} + 228 z^{15} + 494 z^{10} - 228 z^5 + 1)^3}{kz^5(z^{10} - 11 z^5 - 1)^5} = {}
\qquad\qquad\qquad\qquad\qquad\qquad
$$
\begin{equation}
\label{eq:V3M2P5-dodecahedron}
{} = \frac{(z^{30} - 522 z^{25} - 10005 z^{20} - 10005 z^{10} + 522 z^5 + 1)^2}{kz^5(z^{10} - 11 z^5 - 1)^5} + 1.
\end{equation}
Coefficient $k=1728$ is determined from $V^3 = M^2 + kP^5$.

Moreover, we have proved that the Belyi function of a dessin with passport $(3^{20} \,|\, 2^{30} \,|\, 5^{12})$ is equivalent to (\ref{eq:V3M2P5-dodecahedron}). It means that the dodecahedron is the only dessin with its passport.

{\bf Case $s=6$.} In this case $m=12$. The polynomial $V$ 
should be of degree 
$22$, but from (\ref{eq:V3M2P5-V}) we conclude that the coefficient at $z^{22}$ is equal to zero:
$$
\frac{25}{11s^2}(-12 m(m-1) + 11 m^2) \Bigr|_{m=12} = 0
$$
(the same holds for $M$: it should be of degree $33$, but the coefficient at $z^{33}$ is equal to zero). Therefore there is no solutions of (\ref{eq:V3M2P5-Belyi}) satisfying (\ref{eq:V3M2P5-Belyi-deg}) and this completes the proof.
\qed

\medskip

Nevertheless let us check what are the solutions that we found for $s=6$. Plugging (\ref{eq:V3M2P5-P-indeterminate}) into (\ref{eq:V3M2P5-ODE}) and solving sequentially linear equations in $a_i$ we get
$$
a_8 = -\frac{15 a_{10}^2}{44},
\quad
a_7 = -\frac{6 a_9 a_{10}}{55},
\quad
a_6 = -\frac{25 a_{10}^3 + 66 a_9^2}{1210},
$$
$$
a_5 = \frac{3 a_9 a_{10}^2}{1210},
\quad
a_4 = -\frac{3 a_{10} (125 a_{10}^3 + 176 a_9^2)}{106480},
$$
$$
a_3 = -\frac{a_9 (15 a_{10}^3 + 22 a_9^2)}{13310},
\quad
a_2 = \frac{a_{10}^2 (625 a_{10}^3 + 924 a_9^2)}{5856400},
$$
$$
a_1 = \frac{a_9 a_{10} (475 a_{10}^3 + 704 a_9^2)}{64420400},
\quad
a_0 = \frac{3125 a_{10}^6 + 9856 a_9^2 a_{10}^3 + 7744 a_9^4}{2834497600}.
$$
The leading terms of $V$ and $M$ are
$$
V = -\frac{50}{33} a_{10} z^{20} - \frac{50}{11} a_9 z^{19} + \ldots
$$
$$
M = \frac{25}{11} a_9 z^{30} - \frac{2500}{363} a_{10}^2 z^{29} + \ldots
$$
Therefore we have 1-dimensional family of solutions: the two parameters $a_9$ and $a_{10}$ are transformed by the multiplicative substitution $z\to az$ as $(a_9:a_{10})\to (a_9/a^3:a_{10}^2)$, so we can consider the family parametrized by a point of a projective line with the coordinate $(a_9:a_{10})$. One can check that (\ref{eq:V3M2P5-Belyi}) gives us a Belyi function for
$$
k = -\frac{5^4\cdot(2^3\cdot 5^2 \cdot a_{10}^3 + 3^3\cdot 11\cdot a_9^2)}{3^3\cdot 11^3}.
$$
For $a_{10}=0$ we get the dodecahedron dessin with a vertex at the infinite point, while for $a_9=0$ we get the alternative realization of the dodecahedron 
with the middle of an edge placed at the infinite point.

\bigskip

{\bf Generalized Fermat equation and a theorem by Halphen.} The first step in our proof was to reduce the problem of construction of a dessin to the problem of finding complex polynomials satisfying
$$
\frac{V^3}{kP^5} = \frac{M^2}{kP^5} + 1
$$
or
$$
V^3 = M^2 + kP^5,
$$
which is a particular case of the generalized Fermat equation $X^p + Y^q = Z^r$ for $X,Y,Z\in{\mathbb C}[x]$.

The equation attracted a fair amount of interest at the edge of XIX--XX centuries. Both Schwarz and Klein found particular polynomial solutions of the equation. The complete answer was given by Georges Henri Halphen.
In 1880 he wrote a memoir \cite{Halphen-1883} devoted to ordinary linear differential equations and their transformations. The memoir was submitted for a grand prize of the Paris Academy of Sciences and had won the prize.

Although the topic of the memoir is far from diofantine equations, Halphen was naturally led to the generalized Fermat equation in polynomials and obtained the following results:
\begin{enumerate}
\item the equation $X^p + Y^q = Z^r$ has no solution in coprime non-constant polynomials if $1/p+1/q+1/r \le 1$;
\item for any platonian triple $(p,q,r)=(2,2,r)$, $(2,3,3)$, $(2,3,4)$ and $(2,3,5)$ all solutions of the equation has the form
$$
X = aV^{n/p}P(\varphi)),
\quad
Y = bV^{n/q}Q(\varphi)),
\quad
Z = cV^{n/r}R(\varphi)),
$$
where $a,b,c\in {\mathbb C}$ such that $a^n=b^n=c^n$, $\varphi = U/V$ is a rational function, $U,V\in {\mathbb C}[z]$ and $P$, $Q$ and $R$ are the polynomials, in terms of which the Belyi function of the platonian dessin of type $(p,q,r)$ (dihedron, tetrahedron, octahedron, dodecahedron) are expressed.
\end{enumerate}
The former statement is rather easy to prove; one can apply the same argument that Korkine used in his letter to Hermite \cite{Korkine-1880} to prove insolvability of the Fermat equation $X^n+Y^n=Z^n$  in polynomials for $n\ge 3$.

An elementary approach to the problem of finding {\it all} polynomial solutions of the generalized Fermat equation was sketched by Ermakow in 1898 \cite{Ermakow-1898} who challenged young mathematicians to fill in the details. It was done by Welmin (then a student) in 1904 \cite{Welmin-1904}. The dihedron, tetrahedron and octahedron cases turned out to be easy to handle, while the dodecahedron case is very tricky. After 20 pages of lengthy calculations Welmin comes to the answer for $(p,q,r)=(2,3,5)$, but it's completeness is not proved.

Please note that ``elementary'' for Ermakow means that the usage of differentiation is not allowed. In our proof we do not set such artificial constraint and attain conceptually the same result much easier.

See also \cite{ArnaudiesBertin-2001} for a modern exposition of the Halphen's result via Galois theory.

One can apply Halphen's result to deduce our Theorem \ref{th:Non-realizable}. Additional limitations posed by Belyi functions should be used to prove that function $\varphi$ should be of degree 1.

\subsection{$p_6=2$. Barrel}

Specifying the calculations of subsection {\bf 1.2} to the case $p_6=2$, we get
$$
\boxed{\#\text{vertexes}=24,\ \#\text{edges}=36,\ \#\text{faces}=14}
$$

\begin{figure}[ht]
\begin{center}
\includegraphics[scale=1.0]{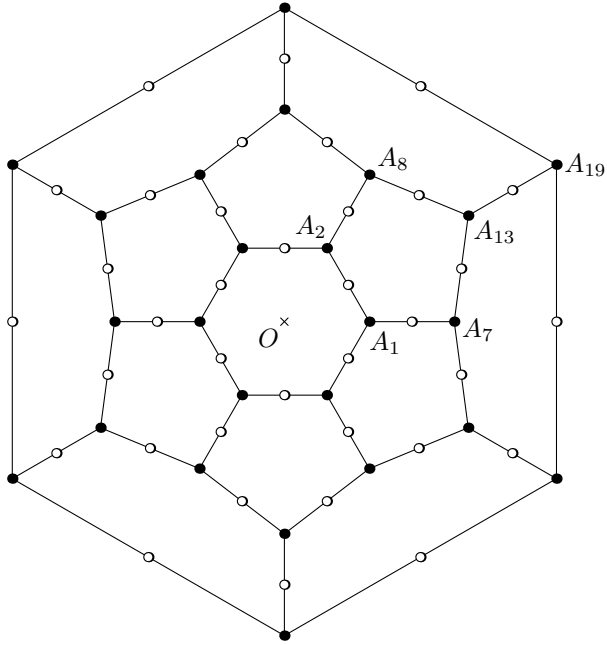}
\end{center}
\caption{Barrel dessin $D_{72}$ (a fulleren with two hexagons).}
\label{fig:barrel}
\end{figure}

It is well known that there exists only one fullerene with these parameters, which is called {\it barrel} and denoted by $C_{24}$. The corresponding dessin $D_{72}$ has 72 edges and is shown in fig.~\ref{fig:barrel}.

The automorphism group of $D_{72}$ is the dihedral group $\mathsf{D}_6$, factorizing by it we get the same $D_{12}$ as in our considerations concerning dodecahdron. It means that we get the Belyi function of $D_{72}$ for free:

\begin{equation}
\label{eq:barrel-Belyi}
\beta_{72}(z) = \beta_{12}(z^6) = \frac{(z^{24} + 228 z^{18} + 494 z^{12} - 228 z^6 + 1)^3}
{1728 z^6 (z^{12} - 11 z^6 - 1)^5}.
\end{equation}

The vertex polynomial of the dessin is
$$
V=z^{24} + 228 z^{18} + 494 z^{12} - 228 z^6 + 1 = \prod_{i=1}^{24}(z-A_i),
$$
where
\begin{equation}
\label{eq:barrel-A}
A_k=\left\{
\begin{array}{ll}
a_1\cdot\mathrm{e}^{\frac{\pi (k-1) \mathrm{i}}{3}}
&\text{ for }k=1,\dots,6;\\
a_7\cdot\mathrm{e}^{\frac{\pi (k-1) \mathrm{i}}{3}}
&\text{ for }k=7,\dots,12;\\
a_{13}\cdot\mathrm{e}^{\frac{\pi  \mathrm{i}}{6}+\frac{\pi (k-1) \mathrm{i}}{3}}
&\text{ for }k=13,\dots,18;\\
a_{19}\cdot\mathrm{e}^{\frac{\pi  \mathrm{i}}{6}+\frac{\pi (k-1) \mathrm{i}}{3}}
&\text{ for }k=19,\dots,24
\end{array}
\right.
\end{equation}
with real positive $a_i$ such that $-a_{19} < -a_{13} < a_1 < a_7$ are the roots of the polynomial $Z^{4} + 228 Z^{3} + 494 Z^{2} - 228 Z + 1$. Finding the roots of this polynomial we get
\begin{equation}
\label{eq:barrel-a}
\begin{array}{lclcl}
a_1 &=& \sqrt[6]{-57-25\sqrt{5}+5\sqrt{255-114\sqrt{5}}} &=& 0.405...,\\[5pt]
a_7 &=& \sqrt[6]{-57+25\sqrt{5}+5\sqrt{255-114\sqrt{5}}} &=& 0.853...,\\[5pt]
a_{13} &=& \sqrt[6]{57-25\sqrt{5}+5\sqrt{255-114\sqrt{5}}} &=& 1.171...,\\[5pt]
a_{19} &=& \sqrt[6]{57+25\sqrt{5}+5\sqrt{255-114\sqrt{5}}} &=& 2.467...
\end{array}
\end{equation}

The formulas (\ref{eq:barrel-A}) and (\ref{eq:barrel-a}) determine the geometry of (our model of) the barrel completely. 
Since we are going to study the pentagon $A_1A_2A_8A_{13}A_7$, we recall
\begin{equation}
\label{eq:barrel-5A}
A_1 = a_1,\quad
A_2=\mathrm{e^{\frac{\pi\mathrm{i}}{3}}}a_1,\quad
A_7=a_7,\quad
A_8=\mathrm{e^{\frac{\pi\mathrm{i}}{3}}}a_7,\quad
A_{13} = \mathrm{e^{\frac{\pi\mathrm{i}}{6}}}a_{13}.
\end{equation}

Please refer to fig.~\ref{fig:barrel} for the positions of the vertexes $A_1$, $A_2$, $A_7$, $A_8$, $A_{13}$. Note that due to the $\mathsf{C}_2$-symmetry that flips the hexagons we have
$$
A_{19}=\frac{\mathrm{e}^{\frac{\pi  \mathrm{i}}{6}}}{A_1},\quad
A_{13}=\frac{\mathrm{e}^{\frac{\pi  \mathrm{i}}{6}}}{A_7}.
$$

\subsection{The geometry  of one pentagonal face}
The positions of the fulleren's vertexes on the standard sphere in $\mathbb{R}^3$ can be calculated by the stereographic projection.

Provide $\mathbb{R}^3$ with the euclidean coordinates $X,Y,Z$ and assume that the standard sphere $\mathbf{S}^2\subset\mathbb{R}^3$ is defined by the equation
$$
X^2+Y^2+Z^2=1.
$$

The above vertexes $A_k$ are understood to lie in the finite part of the projective line
$\mathbf{P}_1(\mathbb{C})$.
The inverse stereographic projection $\mathbb{C}\to\mathbf{S}^2$ is defined (for the real coordinates $z=x+y\mathrm{i}$ on the affine part of the Riemannian sphere $\mathbf{P}_1(\mathbb{C})$) by
$$
\begin{array}{lll}
\displaystyle
X=\frac{2x}{x^2+y^2+1},&
\displaystyle
Y=\frac{2y}{x^2+y^2+1},&
\displaystyle
Z=\frac{x^2+y^2-1}{x^2+y^2+1}.
\end{array}
$$

Using the numerical values given by (\ref{eq:barrel-5A}) one calculates (with the 3-digits precision) the euclidean coordinates of the basic vertexes of the fulleren on the 2-sphere: 
$$
\begin{array}{llcllcllcll}
X(A_1) = 0.696... & Y(A_1) = 0 & Z(A_1) = -0.717...
\\[7pt]
X(A_2) = 0.348... & Y(A_2) = 0.602... & Z(A_2) = -0.717...
\\[7pt]
X(A_7) = 0.987... & Y(A_7) = 0 & Z(A_7) = -0.156...
\\[7pt]
X(A_8) = 0.493... & Y(A_8) = 0.855... & Z(A_8) = -0.156...
\\[7pt]
X(A_{13}) = 0.855... & Y(A_{13}) = 0.493... & Z(A_{13}) = 0.156...
\end{array}
$$

The euclidean lengths of the edges are:
$$|A_1A_7|=.632...,\quad |A_7A_{13}|=.599...,$$
$$|A_{13}A_8|=.599...,\quad |A_8A_2|=.632...,$$
$$|A_2A_1|=.696...$$
(so the axial symmetry of the considered face is confirmed numerically).

It follows from the symmetrical properties of the fulleren that 
{\it the vertexes $A_1,A_2,A_7,A_8$ lie in the same plane}
$$
A_1A_2A_7A_8:\text{ }
\boxed{
0.935... X + 0.540... Y - 0.485... Z = 1
}
$$

The plane containing $A_7,A_8,A_{13}$ is defined by the equation
$$
A_7A_8A_{13}:\text{ }
\boxed{
0.939... X + 0.542... Y - 0.457... Z = 1
}
$$

The angle between these planes equals
$$
\arccos(0.9997...) = 1.36...^\circ,
$$
hence is negligible. 
So {\it the face under consideration is very close to the flat one, but is not actually plane}. We close our study by presenting the plane approximation of the face $A_1A_2A_8A_{13}A_7$.

The angles between its edges are
$$
\angle A_1A_2A_8=\angle A_2A_1A_7
=103.3...^\circ,
$$
$$
\angle A_2A_8A_{13}=\angle A_{13}A_7A_1
=111.2...^\circ,
$$
$$
\angle A_7A_{13}A_8
=110.8...^\circ.
$$

\begin{figure}[h]
\begin{center}
\includegraphics[scale=1.0]{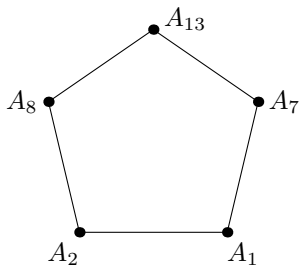}
\end{center}
\caption{Flat approximation of face $A_1 A_7 A_{13} A_8 A_2$ of the barrel.}
\label{fig:barrel-face}
\end{figure}

So the shape of our pentagonal face turns out to be rather close to the regular pentagon. Its flat approximation looks like the one shown in fig.~\ref{fig:barrel-face}.

We are eager to know whether the real barrel is similar to our model.

\section{Discussion}

The most important question is whether the found metrical structure is related to the real geometry of fullerenes. The answer to this question needs some collaboration between mathematicians and chemists. We assume below that the answer is positive (otherwise the existence of {\it two} different ``natural''  metric realizations of the same combinatorial structure will be discovered). Then the following problems occur.

\begin{itemize}
\item Find the {\it physical} representatives of the conformal classes of the general dessins d'enfants, corresponding to fullerenes.
\item Determine in terms of $p_6$ the boundaries of possibilities of modern (super?)computers to find the {\it exact} solutions of the main equation (\ref{eq:2d}).
\item Develop the techniques of the {\it approximate} solutions of the main equation (\ref{eq:2d}).
\item Study the {\it Galois orbits} of fullerenes; it would be exciting to find out that the small orbits (i.e. the ones that can be defined in terms of irrationalities of small degree) coincide with the ones selected by the Nature.
\end{itemize}

\end{document}